\documentclass[12pt]{article}

\usepackage{amsfonts}
\usepackage{float}
\usepackage{flafter}

\leftmargin\leftmargini
\labelwidth\leftmargini\advance\labelwidth-\labelsep

\def\R{\relax\ifmmode I\!\!R\else$I\!\!R$\fi}

\def\Z{\relax\ifmmode Z\!\!\!Z\else$Z\!\!\!Z$\fi}

\def\C{\relax\ifmmode C\!\!\!\!I\else$C\!\!\!\!I$\fi}

\def\K{\relax\ifmmode I\!\!K\else$I\!\!K$\fi}

\def\N{\relax\ifmmode I\!\!N\else$I\!\!N$\fi}

\newcounter{defcounter}[section]
{\vspace{0.1cm}\begin{sloppypar}\noindent\stepcounter{defcounter}{\bfseries Definition
      \thesection.\thedefcounter}}%
{\end{sloppypar}\vspace{0.1cm}}
\newtheorem{corollary}{Corollary}[section]
\newtheorem{lemma}{Lemma}[section]
\newtheorem{theorem}{Theorem}[section]
\newtheorem{proposition}{Proposition}[section]

\newcommand{\proof}{{\noindent\bf Proof. }}

\newcommand{\qed}{\hfill $\square$}

\begin{document}
\thispagestyle{empty}
\begin{center}
{\Large {\bf A construction of singular overlapping asymmetric
self-similar measures}}
\end{center}
\begin{center}
J. Neunh\"auserer\\~\\
Reitstallweg 9, \\ 38640 Goslar, Germany\\
neunchen@aol.com
\end{center}
\begin{abstract}
In \cite{[NE1]} we found a class of overlapping asymmetric
self-similar measures on the real line, which are generically
absolutely continuous with respect to the Lebesgue measure. Here
we construct
exceptional measures in this class being singular. \\
{\bf MSC 2000: Primary 28A12, Secondary 28A78, 28A75, 14H50}
\end{abstract}
\section{Introduction}
In the last century Bernoulli convolutions on the real line, which
may also be described as symmetric self similar measures, were
successfully studied, see \cite{[60]}. Generically these measures
are absolutely continuous with quare integrable density
 and they get singular
for certain algebraic parameters, see \cite{[ER1]} \cite{[ER2]},
\cite{[SO]} and \cite{[PS1]}. In \cite{[NE1]} we began to study
overlapping asymmetric self-similar measures on the reale line,
which generalize Bernoulli convolutions. The parameter domain of
these measures is two dimensional. We extended arguments of Peres
and Solomyak \cite{[PS1]} to show that overlapping asymmetric
self-similar measures are generically (in the sense of Lebesgue
measure on the parameter domain) absolutely continuous with
respect to the Lebesgue measure on the real line (see theorem
2.1). In this article we are concerned with the question, if there
exists any exceptional overlapping asymmetric self-similar
measures being singular. Looking at the symmetric case this is
quiet reasonable. We will prove here that near to the boundary of
the parameter domain there exists exceptional values for which
overlapping asymmetric self-similar measures have a dimension drop
and get hence singular (see theorem 3.2 below). In fact we will
show even more, each point of this boundary is an accumulation
point of parameters for which the dimension of the measures is
less than one. Our strategy to get this result is to proof a
dimension estimate which is effective near to the boundary of the
parameter domain (see theorem 3.1 below). We proof this dimension
estimate using a meanwhile classical approach in ergodic theory.
The dimension of our measures is bounded by the quotient of
entropy and Lyapunov exponent and we are able to estimate the
dimension and calculate the Lyapunov exponent, compare \cite{[PE]}
and especially \cite{[V]}. The assumption of our dimension
estimate is that the parameters values, of the measures we
consider, fulfill a certain type of algebraic equation in two
variables.  We have thus to guarantee the existence of such
solutions in an appropriate domain (see proposition 3.1). To proof
this we use a refinement of construction invented
by Simon and Solomyak in the context of self-similar sets, see \cite{[SS]}.\\
The rest of this paper is organized as follows. In section two we
introduce our notations and basic objects, in section three we
present our results, in section four we proof the dimension
estimate on asymmetric self-similar measures and in section five
we proof the existence of solutions of a certain type of algebraic
equations we used.
\section{Preliminaries}
For $\beta_{1},\beta_{2} \in (0,1)$ consider the linear
contractions
\[ T_{1}x=\beta_{1}x\qquad T_{2}x=\beta_{2}x+1 \]
on the interval $I=[0,1/(1-\beta_{2})]$. For a finite sequence
$s=(s_{1},s_{2},\dots,s_{n}) \in \Sigma_{n}:=\{1,2\}^{n}$ we let
\[T_{s}=T_{s_{1}}\circ T_{s_{2}}\circ \dots \circ T_{s_{n}}\]
Now consider infinite sequences
$s\in\Sigma:=\{1,2\}^{\mathbb{N}}$. We define a map
$\pi:\Sigma\longmapsto I$ by
\[ \pi(s)=\lim_{n\longmapsto \infty }T_{s_{1}}\circ
T_{s_{2}}\circ \dots \circ T_{s_{n}}(I). \]
 For $s\in\Sigma_{n}$
or $s \in\Sigma$ and $k\in\mathbb{N}$ let
\[ 1_{k}(s)=Card\{i|s_{i}=1\mbox{ for }i=1\dots k \} \]
\[ 2_{k}(s)=Card\{i|s_{i}=2\mbox{ for }i=1\dots k \}  \]
By induction we see that for $s\in\Sigma_{n}$ we have
\[ T_{s}x=\beta_{1}^{1_{n}(s)}\beta_{2}^{2_{n}(s)}x+V_{s} \qquad (\star) \]
where
\[ V_{s}=\sum_{k=1}^{n}(s_{k}-1)\beta_{1}^{1_{k-1}(s)}\beta_{2}^{2_{k-1}(s)} .\]
Accordingly the map $\pi$ has the following explicit form,
\[\pi(s)=\sum_{k=1}^{\infty}(s_{k}-1)\beta_{1}^{1_{k-1}(s)}\beta_{2}^{2_{k-1}(s)}. \]
It is well known, see \cite{[HU]} or \cite{[FA]}, that for all
$\beta_{1},\beta_{2}\in(0,1)$ there exists a unique Borel
probability measure $\mu=\mu_{\beta_{1},\beta_{2}}$ on $I$ with
the following self-similarity,
\[ \mu=\frac{T_{1}(\mu)+T_{2}(\mu)}{2}.\]
If $b$ is the equally weighted Bernoulli measure on $\Sigma$ we
get the following description of the equally weighted self-similar
measures $\mu$ using the coding map $\pi$,
\[ \mu=\pi(b)=b\circ \pi^{-1}. \]
If $\beta_{1}=\beta_{2}=\beta\in (0,1)$ the measures
$\mu=\mu_{\beta,\beta}$ are infinite Bernoulli convolutions. As we
mentioned in the introduction  these measures were successfully
studied. In the case $\beta_{2}\not=\beta_{1}$ we call the
measures $\mu=\mu_{\beta_{1},\beta_{2}}$ {\bf asymmetric
self-similar}. We have the following result on the properties of
these measures:
\begin{theorem}
If $\beta_{1},\beta_{2} \in (0,1)$ and $\beta_{1}\beta_{2}<1/4$
than the measure $\mu_{\beta_{1},\beta_{2}}$ is singular with
\[\dim_{H}\mu_{\beta_{1},\beta_{2}}<-2\log2/(\log\beta_{1}+\log\beta_{2})\]
For almost all $\beta_{1},\beta_{2} \in (0,0.649)$ with
$\beta_{1}\beta_{2}\ge 1/4$ the measures
$\mu_{\beta_{1},\beta_{2}}$ are absolutely continuous.
\end{theorem}
This theorem is an immediate consequence of Theorem I in
\cite{[NE1]}. There we worked in a more general setting also
considering measures wich are not equally weighted and studying
the density of these measures. We like to mention that Ngai and
Wang have related results, see \cite{[NW]}. \\If
$\beta_{1}\beta_{2}\ge 1/4$ we call $\mu_{\beta_{1},\beta_{2}}$
{\bf overlapping self-similar}. We do not belief that the bound
$0.649$ in our result on overlapping self-similar measures is
essential, also for some technical reasons we are not able to
remove it. In view of our theorem it is a natural question if
there are any singular overlapping asymmetric self-similar
measures $\mu_{\beta_{1},\beta_{2}}$ at all. Let us first remark
that in the symmetric situation $\beta=\beta_{1}=\beta_{2}$ the
measure $\mu_{\beta,\beta}$ gets singular, if $\beta\in(0.5,1)$ is
the reciprocal of a Pisot number\footnote{An algebraic integer
with conjugates inside the unit circle, see \cite{[ER1]}}. Thus
one may think that there should be some algebraic equations for
$\beta_{1},\beta_{2}\in(0,1)$ with $\beta_{1}\beta_{2}> 1/4$ and
$\beta_{1}\not=\beta_{2}$ such that $\mu_{\beta_{1},\beta_{2}}$
gets singular. We show in this article, that this is in fact true.
\section{Results}
We first present here a new dimension estimate on asymmetric
self-similar measures. By $\dim_{H}\mu$ we denote the Hausdorff
dimension of a measure $\mu$ given by
\[ \dim_{H}\mu=\inf\{\dim_{H} M|\mu(M)=1\} \} \]
where $\dim_{H}M$ is the Hausdorff dimension of a set. See
\cite{[FA]} or \cite{[PE]} for a good introduction to dimension
theory. Now we state our result.
\begin{theorem}
Let $\beta_{1},\beta_{2}\in(0,1)$ and $s,t\in\Sigma_{n}$ with
$s\not= t$.  If $T_{s}=T_{t}$ then
\[ \dim_{H}\mu_{\beta_{1},\beta_{2}}<-\frac{2
\log(2^{n}-1)}{n(\log\beta_{1}+\log\beta_{2})}.\]
\end{theorem}
This dimension estimate has the following obvious corollary on the
singularity of asymmetric self-similar measures.
\begin{corollary}
Let $\beta_{1},\beta_{2}\in(0,1)$ and $s,t\in\Sigma_{n}$ with
$s\not= t$. If $T_{s}=T_{t}$ and
\[ \frac{1}{(\sqrt[n]{2^{n}-1})^{2}} > \beta_{1}\beta_{2}  \]
then $\mu_{\beta_{1},\beta_{2}}$ is singular with
$\dim_{H}\mu_{\beta_{1},\beta_{2}}<1$.
\end{corollary}
We remain here with the problem if there are any
$(\beta_{1},\beta_{2})$ in the domain
\[ D_{n} =\{ (\beta_{1},\beta_{2})\in (0,1)|\frac{1}{(\sqrt[n]{2^{n}-1})^{2}} >
\beta_{1}\beta_{2}> \frac{1}{4}\}\] such that $T_{s}=T_{t}$ for
$s,t\in\Sigma_{n}$ with $s\not= t$. By equation $(\star)$ on the
last side this is equivalent to find solutions of an algebraic
equations $V_{s}=V_{t}$ with $s\not= t$ and $1_{n}(s)=1_{n}(t)$ in
the domain $D_{n}$. We will adopt a construction developed by
Simon and Solomayak \cite{[SS]} to prove the following
proposition.
\begin{proposition}
For all $\beta_{2} \in(1/4,1/2)$ there exists constants $c>0$,
$\lambda\in(0,1)$, a sequence $n_{k}\longmapsto\infty$, an
\[\beta_{1}\in
(\frac{1}{4\beta_{2}},\frac{1}{4\beta_{2}}+c\lambda^{n_{k}}) \]
and $s,t\in\Sigma_{n_{k}}$ with $s\not= t$ such that
$T_{s}=T_{t}$.
\end{proposition}
By this proposition and corollary 3.1 we get the main result of
the articel
\begin{theorem}
For all $\beta_{2} \in(1/4,1/2)$ and $\epsilon>0$ sufficient small
there is an
\[ \beta_{1}\in
(\frac{1}{4\beta_{2}},\frac{1}{4\beta_{2}}+\epsilon)\]  such that
the measure $\mu_{\beta_{1},\beta_{2}}$ is singular with
$\dim_{H}\mu_{\beta_{1},\beta_{2}}<1$.
\end{theorem}
\proof Fix $\beta_{2}\in (1/4,1/2)$ and $c,\lambda$ from
proposition 3.1. By the rule of Le Hospital the convergence of
$(\sqrt[n]{2^{n}-1})^{2}$ to $4$ is subexponential hence there
exists $n_{0}\in\mathbb{N}$ such that
\[ \forall n\ge n_{0}:\frac{1}{4}+c\beta_{2}\lambda^{n} < \frac{1}{(\sqrt[n]{2^{n}-1})^{2}}\]
Choose in proposition 3.1 an $n_{k}>n_{0}$ and fix
$\epsilon_{0}=c\lambda^{n_{k}}$. Now by proposition 3.1 for
$\epsilon<\epsilon_{0}$ there is an $n_{k}>n_{0}$ with
$c\lambda^{n_{k}}<\epsilon$ and there exists
\[\beta_{1}\in
(\frac{1}{4\beta_{2}},\frac{1}{4\beta_{2}}+\epsilon)\] with
$T_{s}=T_{t}$ for $s,t\in\Sigma_{n_{k}}$. In addition we have
\[ 1/4<\beta_{1}\beta_{2}<1/4+\beta_{2}c\lambda^{n_{k}}<\frac{1}{(\sqrt[n_{k}]{2^{n_{k}}-1})^{2}} \]
Hence $(\beta_{1},\beta_{2})\in D_{n_{k}}$ and by corollary 3.1 we
get $\dim_{H}\mu_{\beta_{1},\beta_{2}}<1$ implying singularity of
the measure $\mu_{\beta_{1},\beta_{2}}$ .\qed ~\\~\\
 We may
paraphrase our result as follows. In the parameter domain where
overlapping asymmetric self similar measures are generically
absolutely continues we find exceptional parameter values,
sufficient near to any point of the boundary, for which these
measures are singular. Each point of this boundary is an
accumulation point of such exceptional values. Of course our
result leaves it open if there are any exceptional values near to
points with a given distance from the boundary. Our dimension
estimate in theorem 3.1 is not strong enough to answer this
question. We conjecture from an analogy with Bernoulli
convolutions that this estimate is not sharp in general (compare
with \cite{[NE2]}).
\section{A dimension estimate for asymmetric self-similar measures}
Our aim in this section is to proof theorem 3.1 using a classical
dimension estimate by entropy and Lyapunov exponent. To this end
we first introduce the Lyapunov exponent $\Xi(\mu)$ of the
self-similar measure $\mu=\mu_{\beta_{1},\beta_{2}}$. For
$s\in\Sigma$ let
\[ \Xi_{s}(\mu)=\lim_{n\longmapsto\infty}\frac{1}{n}\log||T_{s_{1}}\circ T_{s_{2}}\circ \dots \circ
T_{s_{n}}||=\lim_{n\longmapsto\infty}\frac{1}{n}\log\beta_{1}^{1_{n}(s)}\beta_{2}^{2_{n}(s)}\]
The following lemma is just an application of Birkhoff`s ergodic
theorem, see for instance \cite{[KH]}.
\begin{lemma}
For allmost all $s\in\Sigma$ with respect to $b$ we have
\[ \Xi_{s}(\mu)=-\frac{1}{2}(\log\beta_{1}+\log \beta_{2})=:\Xi(\mu) \]
\end{lemma}
\proof Just apply the ergodic theorem to the function
$f:\Sigma\longmapsto \mathbb{R}$
\[ f(s)=\{\begin{array}{cc} \log\beta_{1}\quad\mbox{ for }\quad s_{1}=1 \\ \log\beta_{2}\quad\mbox{ for }\quad s_{1}=2  \end{array}   \]
\qed~\\
Now we introduce the entropy $h(\mu)$ of the self-similar measure
$\mu$. For $m\in\mathbb{N}$ define a Partition ${\bf P}_{n}$ of
$\Sigma$ by the following equivalence relation
\[ s\sim t :\Leftrightarrow T_{s_{1}}\circ T_{s_{2}}\circ \dots \circ
T_{s_{m}}=T_{t_{1}}\circ T_{t_{2}}\circ \dots \circ T_{t_{m}}.\]
The entropy of the partition is given by
\[H({\bf P}_{m})=-\sum_{P\in{\bf P}_{n}}b(P)\log b(P). \]
The sequence $H({\bf P}_{n})$ is subadditive hence the limit
\[ h(\mu):=\lim_{m\longmapsto\infty }\frac{H({\bf P}_{m})}{m} \]
exists. In the following lemma we state a simple estimate on the
entropy $h(\mu)$.
\begin{lemma}
If $T_{s}=T_{t}$ holds for $s,t\in\Sigma_{n}$ with $s\not= t$,
then we have
\[ h(\mu)\le \frac{1}{n}\log(2^{n}-1). \]
\end{lemma}
\proof It is well knowen that
\[ H({\bf P}_{m})\le \log(Card({\bf P}_{m})) \]
see  for instance \cite{[KH]}. By the definition of ${\bf P}_{m}$
and our assumption we have
\[ Card({\bf P}_{nm})\le (2^{n}-1)^{m}\]
Hence
\[ h(\mu)=\lim_{m\longmapsto\infty}\frac{H({\bf
P}_{nm})}{nm}\le\lim_{m\longmapsto\infty}\frac{\log(Card({\bf
P}_{nm}))}{nm} =\frac{1}{n}\log(2^{n}-1)\] \qed \\~\\
In the dimensiontheory it has be shown that for a broad classes of
measures the Hausdorff dimension is bounded from above by the
quotient of entropy and Lyapunov exponent, see \cite{[PE]}. For
symmetric self-similar measures this result is contained in
\cite{[LA]} and for general classes measures constructed by
iterated function systems it is proofed in \cite{[V]}. The
following proposition is just an application of theorem 2.2
in\cite{[V]} to the iterated function system $(T_{1},T_{2})$.
\begin{proposition} With the notations from above
\[ \dim_{H}\mu \le \frac{h(\mu)}{\Xi(\mu)}\]
\end{proposition}
Theorem 3.1 is now obvious from the results of this section.
\section{Existence of solutions of certain algebraic equations in two
variables} In this section we will proof proposition 3.1. We have
to find sequences $s,t\in\Sigma_{n}$ with $s\not= t$ such that
$T_{s}=T_{t}$ holds for some $\beta_{1},\beta_{2}\in(0,1)$. We
will use here a slight modification of a construction we found in
the proof of proposition
3.4. of Simon and Solomyaks work on self-similar sets \cite{[SS]}.\\
Let $s\in\Sigma_{n}$ be the sequence given by
$12^{N_{1}}1^{M_{1}}\dots 1^{M_{k-1}}2^{N_{k}}1^{M_{k}}$ where we
use the exponential here to describe the number of repetitions of
an entry. Let ${\mathcal N}={\mathcal N}_{k}=N_{1}+\dots N_{k}$,
${\mathcal M}={\mathcal M}_{k}=1+M_{1}+\dots+M_{k}$ and
$n=n_{k}={\mathcal M}_{k}+{\mathcal N}_{k}$. Furthermore let
$t\in\Sigma_{n}$ be given by $21^{{\mathcal M}+1}2^{{\mathcal
N}-1}$. By the definition of $T_{s}$ and $T_{t}$ we get:
\[ T_{s}x=\beta_{1}^{{\mathcal N}}\beta_{2}^{{\mathcal M}}x+V_{s}(\beta_{1},\beta_{2})\quad\mbox{ and }\quad T_{t}x=\beta_{1}^{{\mathcal N}}\beta_{2}^{{\mathcal M}}x+V_{t}(\beta_{1},\beta_{2}) \]
with
\[V_{s}(\beta_{1},\beta_{2})=\sum_{j=0}^{k-1}(\beta_{2}^{N_{1}+\dots+N_{j}}\beta_{1}^{1+M_{1}+\dots+M_{j}}\sum_{i=0}^{N_{j+1}-1}\beta_{2}^{i})  \]
\[V_{t}(\beta_{1},\beta_{2})=1+\beta_{1}^{{\mathcal M}}\sum_{i=1}^{{\mathcal N}-1}\beta_{2}^{i} \]
Obviously we have $T_{s}=T_{t}$ if and only if the algebraic
equation \[ V_{s}(\beta_{1},\beta_{2})=V_{t}(\beta_{1},\beta_{2})
\] holds. We first construct by recursion a sequences of natural
numbers $M_{l}$ and $N_{l}$ such that $V_{s}(\beta_{1},\beta_{2})$
is arbitrary close to $V_{t}(\beta_{1},\beta_{2})$. Let
\[ N_{1}=\max\{N\ge
1|y_{1}:=\beta_{1}\sum_{i=0}^{N-1}\beta_{2}^{i}<1\}\]
\[ M_{1}= \min\{M\ge 1|y_{1}+\beta_{2}^{N_{1}}\beta_{1}^{M+1}<1\} \]
If $N_{l-1}$ $M_{l-1}$ an $y_{l-1}$ are constructed continue with
\[ N_{l}=\max\{N\ge 1|y_{l}=y_{l-1}+\beta_{2}^{N_{1}+\dots+N_{l-1}}\beta_{1}^{1+M_{1}+\dots+M_{l-1}}\sum_{i=0}^{N-1}\beta_{2}^{i}<1\}\]
\[ M_{l}=\min\{M\ge 1|y_{l}+\beta_{2}^{N_{1}+\dots+N_{l}}\beta_{1}^{1+M_{1}+\dots+M_{l-1}+M}<1\}\]
By this construction we achieve:
\begin{lemma}
With the notion from above we have for all $\beta_{1},\beta_{2}\in
(0,1)$ with $\beta_{1}+\beta_{2}>1$ and $k\ge 2$
\[ N_{k}\le
\max\{1,\log_{\beta_{2}}(\frac{\beta_{1}+\beta_{2}-1}{\beta_{1}})\}=:B
\] and
\[ V_{t}(\beta_{1},\beta_{2})-V_{s}(\beta_{1},\beta_{2})\le (\frac{1}{\beta_{1}}+\frac{\beta_{2}}{1-\beta_{2}})(\sqrt[1+B]{\beta_{1}})^{n_{k}} \]
\[ \frac{\partial{(V_{t}(x,\beta_{2})-V_{s}(x,\beta_{2}))}}{\partial x} \le
\frac{\beta_{2}}{1-\beta_{2}}{\mathcal M}_{k}x^{{\mathcal
M}_{k}-1}-1\]
\end{lemma}
\proof By the definition of $M_{k-1}$ we have
\[ y_{k-1}+\beta_{2}^{{\mathcal N}_{k-1}}\beta_{1}^{{\mathcal
M}_{k-1}-1}\ge 1 \] and by the definition of $N_{k}$ we have
\[y_{k-1}+\beta_{2}^{{\mathcal N}_{k-1}}\beta_{1}^{{\mathcal
M}_{k-1}}\frac{1-\beta_{2}^{N_{k}}}{1-\beta_{2}}<1.\] Hence
\[\beta_{2}^{{\mathcal N}_{k-1}}\beta_{1}^{{\mathcal
M}_{k-1}}\frac{1-\beta_{2}^{N_{k}}}{1-\beta_{2}}< 1-y_{k-1}\le
\beta_{2}^{{\mathcal N}_{k-1}}\beta_{1}^{{\mathcal M}_{k-1}-1}\]
and thus
\[ \beta_{2}^{N_{k}}\ge \frac{\beta_{1}+\beta_{2}-1}{\beta_{1}}
\]
given the bound on $N_{k}$ state in our lemma. By this we
immediately get
\[ {\mathcal N}_{k}\le Bk\le B{\mathcal M}_{k}\]
and
\[ n_{k}= {\mathcal N}_{k}+ {\mathcal M}_{k}\le (1+B){\mathcal M}_{k}.\]
With $V_{t}(\beta_{1},\beta_{2})<1+\beta_{1}^{{\mathcal
M}_{k}}(\beta_{2}/(1-\beta_{2}))$ and
$V_{s}(\beta_{1},\beta_{2})\ge 1-\beta_{2}^{{\mathcal
N}_{k}}\beta_{1}^{{\mathcal M}_{k}-1}$ we now estimate
\[ V_{t}(\beta_{1},\beta_{2})-V_{s}(\beta_{1},\beta_{2})\le
\beta_{2}^{{\mathcal N}_{k}}\beta_{1}^{{\mathcal M}_{k}-1}
+\beta_{1}^{{\mathcal M}_{k}}\frac{\beta_{2}}{1-\beta_{2}}\le
\beta_{1}^{{\mathcal
M}_{k}}(\frac{1}{\beta_{1}}+\frac{\beta_{2}}{1-\beta_{2}})
\]
\[ \le(\frac{1}{\beta_{1}}+\frac{\beta_{2}}{1-\beta_{2}})(\sqrt[1+B]{\beta_{1}})^{(1+B){\mathcal
M}_{k}}
\le(\frac{1}{\beta_{1}}+\frac{\beta_{2}}{1-\beta_{2}})(\sqrt[1+B]{\beta_{1}})^{n_{k}}.
\]
Calculating derivatives gives
\[V_{s}(x,\beta_{2})=\sum_{j=0}^{k-1}(\beta_{2}^{{\mathcal N_{j}}}x^{{\mathcal M_{j}}}\sum_{i=0}^{N_{j+1}-1}\beta_{2}^{i})  \]
\[V_{t}(x,\beta_{2})=1+x^{{\mathcal M}}\sum_{i=1}^{{\mathcal N}-1}\beta_{2}^{i}. \]
Hence
\[\frac{ \partial V_{s}(x,\beta_{2})}{\partial x}\ge 1\mbox{  and
} \frac{\partial V_{t}(x,\beta_{2})}{\partial
x}<\frac{\beta_{2}}{1-\beta_{2}}{\mathcal M}_{k}x^{{\mathcal
M}_{k}-1}
\]
given the estimate on the derivative. \qed ~\\~\\ Now we are
prepared to prove a proposition on the solution of algebraic
equations $V_{s}=V_{t}$.
\begin{proposition}
For all $\beta_{1},\beta_{2}\in (0,1)$ with $\beta_{1}+\beta_{2}>
1$ there exits $n_{k}\longmapsto \infty$ and
$s,t\in\Sigma_{n_{k}}$ such that the equation
\[V_{s}(x,\beta_{2})=V_{t}(x,\beta_{2}) \]
has a solution
\[x\in (\beta_{1},\beta_{1}+2(\frac{1}{\beta_{1}}+\frac{\beta_{2}}{1-\beta_{2}})(\sqrt[1+B]{\beta_{1}})^{n_{k}}). \]
\end{proposition}
\proof Let \[ f(x)=V_{s}(x,\beta_{2})-V_{t}(x,\beta_{2}) \] and
\[ C_{k}=
(\frac{1}{\beta_{1}}+\frac{\beta_{2}}{1-\beta_{2}})(\sqrt[1+B]{\beta_{1}})^{n_{k}}.
 \]
By our construction of the sequences $s,t\in\Sigma_{n_{k}}$ and
lemma 5.1 $f(\beta_{1})\le C_{k}$. Moreover by lemma 5.1
$f'(x)<-1/2$ for all $x\in [\beta_{1},\beta_{1}+2C_{k}]$ if $k$ is
large enough. By elementary calculus we have $f(x)=0$ for
$x\in(\beta_{1},\beta_{1}+2C_{k})$ proving our proposition.
\qed~\\~\\
With the notions of this section proposition 3.1 is just a
corollary to the last proposition.~\\~\\
 {\bf Proof of proposition 3.1.}
Let $\beta_{2}\in(1/4,1/2)$ and $\beta_{1}=1/(4\beta_{2})$.
Obviously $\beta_{1}+\beta_{2}>1$. Note that $T_{s}=T_{t}$ if
$V_{s}=V_{t}$. Hence proposition 5.1 directly implies Proposition
3.1 with
\[c=2(\frac{1}{\beta_{1}}+\frac{\beta_{2}}{1-\beta_{2}})\mbox{ and
 }\lambda=\sqrt[1+B]{\beta_{1}}.\]
\qed

\end{document}